\documentclass[11pt,reqno]{amsart}

\usepackage[T1]{fontenc}
\usepackage{graphicx}
\usepackage{cite,enumerate}

\usepackage{color}
\definecolor{MyLinkColor}{rgb}{0,0,0.4}

\newcommand{\id}{\mathop{\rm id}\nolimits}
\newcommand{\tr}{\mathop{\rm tr}\nolimits}
\newcommand{\p}{\partial}

\newcommand{\f}{\varphi}

\newcommand{\h}{\rho}

\newcommand{\ov}{\overline}
\newcommand{\wh}{\widehat}
\newcommand{\A}{\mathcal{A}}
\newcommand{\B}{\mathcal{B}}
\newcommand{\kH}{\mathcal{H}}
\newcommand{\0}{\Omega}

\newcommand{\Ro}{{\rm Re}}

\newcommand{\E}{\mathcal{E}}
\newcommand{\V}{\mathcal{V}}

\newcommand{\kL}{\mathcal{L}}

\newcommand{\s}{\mathbb{S}}

\newcommand{\R}{\mathbb{R}}

\newcommand{\W}{\mathcal{W}}

\newcommand{\C}{\mathbb{C}}
\newcommand{\x}{\mathcal{T}}
\newcommand{\mS}{\mathcal{S}}
\newcommand{\N}{\mathbb{N}}

\newcommand{\Z}{\mathbb{Z}}

\newtheorem{thm}{Theorem}[section]
\newtheorem{prop}[thm]{Proposition}
\newtheorem{lem}[thm]{Lemma}

\newtheorem{obs}[thm]{Observation}
\theoremstyle{remark}

\numberwithin{equation}{section}

\title[Well-posedness and stability properties]{Well-posedness and stability analysis for a moving boundary problem modelling the growth of nonnecrotic tumors}


\usepackage[colorlinks=true,linkcolor=MyLinkColor,citecolor=MyLinkColor]{hyperref}

\author[J. Escher]{Joachim Escher}
\address{Institut f{\"u}r Angewandte Mathematik, Leibniz Universit{\"a}t Hannover, Welfengarten~1, 30167 Hannover, Germany. }
\email{escher@ifam.uni-hannover.de}

\author[A.-V. Matioc]{Anca-Voichita Matioc}
\address{Institut f{\"u}r Angewandte Mathematik, Leibniz Universit{\"a}t Hannover, Welfengarten~1, 30167 Hannover, Germany. }
\email{matioca@ifam.uni-hannover.de}

\usepackage[colorlinks=true,linkcolor=MyLinkColor,citecolor=MyLinkColor]{hyperref} 

\begin{document}

\begin{abstract}
We study a moving boundary problem describing the growth of nonnecrotic
tumors in different regimes of vascularisation.
This model  
consists of two decoupled Dirichlet  problem, one for the rate at which nutrient is added to the tumor domain and one for the pressure 
inside the tumor. 
These variables are coupled by a relation which describes the dynamic of the boundary.
By 
re-expressing the problem  as an abstract evolution equation, we prove 
local well-posedness in the small H\"older spaces context. 
Further on,  we use the principle of  linearised stability to characterise the stability properties of the unique radially symmetric equilibrium of the problem. 
\end{abstract}

\maketitle

\section{Introduction}
The study of tumor growth models  is a very current topic in mathematics.
During the last four decades 
an increasing number of mathematical models have been proposed to describe the growth of solid tumors (see \cite {BLM, Cui, CE, FR} and the literature therein).
There is a three level approach in modeling the complex phenomena influencing and describing the processes inside a tumor.
 Models at {\it sub-cellular level} take into consideration that the evolution of a cell is determined by the genes in its nucleus, at {\it cellular level} they model cell-cell interaction and at {\it macroscopic level}, when the tumor is considered to consist of three zones: an external proliferating zone near high concentration of nutrient, an intermediate layer and an internal zone consisting of necrotic cells only. Very often  models combine aspects from different scales.
There are also, a large variety of different types of models: {\it biological models,} consisting of coupled ODE systems where the variables correspond to some biological properties of an entire population; {\it mechanical models} yield to determine the cell movement based on physical forces; {\it the discrete models } handle single-cell scale phenomena and the effects are then examined at macroscopic scale and {\it moving boundary models }  when the macroscopic description of biological tissues is obtained from continuum   mechanics or microscopic description at cellular level. 

In this paper we deal with a moving boundary problem, which is obtained by combining aspects from the cellular and macroscopic scale, and possesses also characteristics of the mechanical model (Darcy's law). Cristini et al. obtained in \cite{VCris}, using algebraic manipulations, 
a new mathematical formulation of an existing model (see \cite{BC, FR, Gr}), which describes the  evolution of nonnecrotic tumors in all 
regimes of vascularisation.
This new formulation has the advantage of considering different  intrinsic-time and length-scales related to  the evolution of the tumor
 and, by incorporating them in the modeling, provides a model describing both vascular and avascular tumor:
\begin{equation}\label{eq:problem}
\left \{
\begin{array}{rlllll}
\Delta \psi &=& f(\psi )  &\text{in} \ \Omega(t) ,&t\geq0, \\[1ex]
\Delta p &=& 0 & \text{in}\ \Omega(t),&t\geq0,  \\[1ex]
\psi &=& 1 & \text{on}\ \partial \Omega(t),&t\geq0,  \\[1ex]
p&=& \kappa_{\p\0(t)}- AG \displaystyle\frac{ |x|^2}{4} &\text{on} \ \partial
 \Omega(t),&t\geq0,\\[1ex]
G\displaystyle\frac{\p \psi}{\p n} -\displaystyle\frac{\p p}{\p n} -AG \displaystyle\frac{n\cdot x}{2} &=&V(t)
 & \text{on}\ \partial \Omega(t),&t>0,  \\[1ex]
\Omega(0)&=&\Omega_0.&&
\end{array}
\right.
\end{equation}
Hereby  $\0_0$ is the initial state of the tumor,
$V$ is the normal velocity of the tumor boundary, $\kappa_{\p\0(t)}$  the curvature of $\p\0(t),$ and the constants $A$ and $G$ have biological meaning, namely $G$ is the rate of mitosis (cell proliferation) and $A$ describes the balance between the rate of mitosis and apoptosis (naturally cell death).
The function $f\in C^\infty([0,\infty))$ has the following properties
\begin{equation}\label{eq:conditions}
f(0)=0 \qquad\text{and}\qquad f'(\psi)>0 \quad\text{for}\quad \psi\geq
 0.
 \end{equation}
 The tumor domain $\0(t)$ is an unknown of the problem and, together with   
 the rate $\psi$ at which nutrient is added to the tumor domain $\0(t)$ and  the pressure $p$ inside the tumor, 
 is to be determined.
 
 Three different  regimes of vascularisation are introduced by the constants $A$ and $G:$ if $G\geq0$ and $A>0 $ the tumor is low vascularised,
 $G\geq0$ and $A\leq0$ correspond to the moderate vascularised case, and if $G<0$ the tumor is highly vascularised.
  
In \cite{VCris} the special case  $f=\id_{[0,\infty)}$ is analysed numerically. 
Moreover, in this situation, the first equation of the system is linear, and 
if the tumor domain is a sphere or an infinite cylinder, then the solution is known through an explicit formula.
The  radially symmetric case when tumors are circles is considered in \cite{EM}, where we show that if $A\in(0,f(1)),$ then there exists a
unique radially symmetric stationary solution $D(0,R_A)$ of \eqref{eq:problem}.
The radius of the stationary solution depends only on $A$ and this circular steady-state is exponentially stable under radially symmetric perturbations in the avascular case  when $G>0,$ and unstable in the high 
vascularisation regime $G<0$, result established for  $f=\id_{[0,\infty)}$ also in \cite{VCris}.
The analysis in \cite{EM} will serve us in the present paper as an ancillary tool when proving the local well-posedness
of problem \eqref{eq:problem} and when studying the stability properties of $D(0,R_A)$.

The model,  presented in \cite{BC, FR, Gr}, has been studied extensively by different authors, see e.g. \cite{BF05, Cui,  CE1, CE, FR, FR01} and the references therein. In particular, it is shown in these papers that if certain parameters belong to an appropriate range, then the mathematical formulation possesses 
a unique radially symmetric solution, result matching perfectly with \cite{EM}.
 Moreover, the stability properties of this solution under general perturbations, as well as bifurcation phenomena are studied.
In contrast, for the  model presented in \cite{VCris, EM}, and which we consider herein, not many analytic results are available.
We prove that also this model is locally well-posed in time, meaning that for appropriate smooth initial data $\0_0$, there exists a unique solution of 
\eqref{eq:problem}, cf. Theorem \ref{Thm1.1}. 
Though in the radially symmetric case the steady-state solution $D(0,R_A)$ is exponentially stable if $G>0,$
we show in Theorem \ref{P:G*G}, by considering arbitrary initial data,
 that this solution is unstable also in the low vascularised case, provided $G$ lies above a well-defined constant $G_*.$
This result matches the case $\gamma<\gamma_* $ in \cite[Theorem 1.2]{CE1}, since $G$ is inversely   proportional to $\gamma.$
The situation when $G\in(0,G_*)$ is still an open problem.   
If $G=0$ the problem is equivalent to the  Hele-Shaw problem studied in \cite{EM4} and 
the exponential stability result stated by \cite[Theorem 4.3]{EM4} holds true.
As a new property we establish in Theorem \ref{T:conv} exponential convergence of $D(0,R_A)$ for every $G>0$  and initial data in a certain  class which  depends on $G$.
 
The outline of the paper is as follows: we introduce in the second section a parametrisation for the unknown tumor domain which permites us to present  the main results Theorems \ref{Thm1.1}-\ref{T:conv}.
Section 3 is dedicated to the proof of Theorem \ref{Thm1.1}, and the stability results stated in Theorems \ref{P:G*G} and \ref{T:conv} are
proved in Section 4.

\section {The main results}
Let $R>0$ be fixed for the remainder of this section. 
Our goal is to show that if the tumor is initially close to $D(0,R),$ then problem \eqref{eq:problem}
 possesses a unique classical H\"older solution.
 To this scope let $h^{r}(\s)$, $r\geq 0$, 
denote the closure of the smooth functions $C^\infty(\s)$
in the H\"older space $C^{r}(\s)$.
Hereby, $\s$ stands for the unit circle and we identify functions on $\s$ with $2\pi$-periodic functions on $\R.$
The  small H\"older spaces $h^{r}(\s)$ have the nice property that the embedding 
$h^{r}(\s)$ is densely and compactly in $h^{s}(\s)$  
for all $ 0\leq s<r.$
We fix $\alpha\in(0,1)$ and we shall use functions $\h\in\V,$ 
whereby
\begin{equation*}
\V:=\{ \rho\in h^{4+\alpha}(\s)\,:\, \|\rho\|_{C(\s)}<1/4 \},
\end{equation*}
to parametrise the boundary of the tumor domain.
Obviously, $\V$ is an open neighbourhood  of the zero function    in $ h^{4+\alpha}(\s)$.
Given $\rho\in\V,$ we define the $C^{4+\alpha}$-perturbation of the  circle centred in $0$ with radius $R$
\[
\Gamma_\rho:=\left\{x\in\R^2\,:\,|x|=R\left(1+\rho\left(x/|x|\right)\right)\right\}
=\left\{R\left(1+\rho(x)\right)x\,:\,x\in\s\right\}.
\]
The simply connected component of $\R^2$ which is bounded by the curve $\Gamma_\rho$ is the set
\[
\0_\rho:=\left\{x\in\R^2\, :\, |x|<R\left(1+\rho\left(x/|x|\right)\right)\right\}\cup\{0\},
\]
with boundary   $\p\0_\rho=\Gamma_\rho.$
Given $x\in\Gamma_\rho,$ the real number $\rho(x/|x|)$ is the ratio of the signed distance from $x$ to the circle $R\cdot \s$ and $R$ (see Figure 1). 
\begin{figure}\label{F:domeniu2}
$$\includegraphics[width=0.55\linewidth]{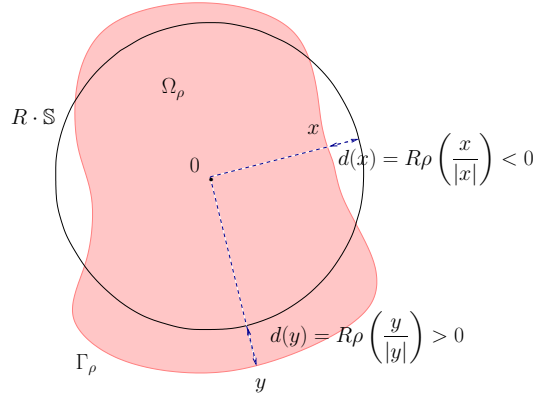}$$
\caption{Parametrisation of the tumor domain}
 \end{figure}
It is suitable to represent $\Gamma_\rho $ as the $0-$level set of an appropriate function.
For this, let $N_\rho: A(3R/4,5R/4)\to\R$ be the function defined by
\[
N_\rho(x)=|x|-R-R\rho(x/|x|), \quad x\in A(3R/4, 5R/4),
\]
where $A(3R/4,5R/4)$ is the annulus centred in $0$ with radii $3R/4$ and $5R/4$
$A(3R/4,5R/4):=\{x\in\R^2:\, 3R/4<|x|<5R/4\}.$
Obviously $A(3R/4,5R/4)$ is an open neighbourhood of $\Gamma_\rho $ and
 $\Gamma_\rho=N^{-1}_\rho(0)$.
Let  $\nu_\rho$ denote  the outward normal at $\Gamma_\rho.$ 
Since $\Gamma_\rho$ is the $0-$level set of $N_\rho$,  the gradient $\nabla N_\rho$ and $\nu_\rho$ must be collinear vectors.
Moreover, $N_\h$ is positive on the complement of   $\ov\0_\rho $, hence $N_\rho(x+\lambda\nu_\rho(x))>0$ for all $x\in\Gamma_\rho $ and $\lambda>0.$
Differentiating this relation with respect to $\lambda$ at $\lambda=0$ yields that $\nabla N_\rho\cdot \nu_\rho>0,$ hence
\[
\nu_\rho=\frac{\nabla N_\rho}{|\nabla N_\rho|}.
\]
To incorporate time let $T>0.$
Presuppose that the function $\rho\in C([0,T],\V)\cap C^{1}([0,T],h^{1+\alpha}(\s))$ describes the evolution of the tumor,  which at time 
$t=0$ is located at $\0(0)=\0_{\rho(0)}.$
The normal velocity  $V(t)$ of the moving boundary $\Gamma_{\rho(t)}$ is then given by the expression
\[
V(t)=-\frac{\p_tN_\rho}{|\nabla N_\rho|}.
\]
This relation follows from the standard assumption that the interface moves along with the tumor and from  relation $\Gamma_\rho=N_\rho^{-1}(0).$
With this notation \eqref{eq:problem} is equivalent  to the following system of equations
\begin{equation}\label{3}
\left \{
\begin{array}{rlcl}
\Delta \psi &=& f(\psi )  &\text{in} \ \Omega _{\rho(t)},  \\[1ex]
\Delta p &=& 0 & \text{in}\ \Omega _{\rho(t)},  \\[1ex]
\psi &=& 1 & \text{on}\ \Gamma _{\rho(t)},  \\[1ex]
p&=& \kappa_{\Gamma_{\rho(t)}}- AG \displaystyle\frac{ |x|^2}{4} &\text{on} \ \Gamma _{\rho(t)}, \\[1ex]
\partial_t N_{\rho}&=& - \left<G\nabla \psi -\nabla p- AG\displaystyle \frac{ x}{2}, \nabla N_{\rho}\right>  & \text{on}\ \Gamma _{\rho(t)}, \\[1ex]
\rho (0) &=& \rho_0  & \text{on}\ \s,
\end{array}
\right.
\end{equation}
for all $t\in [0,T]$.
A triple $(\h, \psi, p)$ is called a {\em classical H\"older solution} of \eqref{3} if
\begin{align*}
&\text{$\rho \in C([0,T], \V)\cap C^1([0,T], h^{1+\alpha}(\s)),$}\\[1ex]
& \text{$\psi(\cdot, t)p(\cdot, t)\in buc^{2+\alpha}(\Omega_{\rho(t)})$ for all $t\in[0,T],$}
\end{align*}
and $(\h, \psi, p)$ solves the system \eqref{3} pointwise.
Given $\h\in\V,$    $ buc^{k+\alpha}(U)$ stands for   the closure of $ BU\!C^\infty(U)$
in $ BU\!C^{k+\alpha}(U).$
Of major interest  is to determine the mapping $\rho$  which describes the evolution of the tumor.
The functions $\psi$ and $p$ can be  then determined as solutions of Dirichlet problems, cf. Lemmas \ref{Lem2.1},  \ref{L.2.3}, and \ref{L.2.4}. 
This is the reason why we shall also refer only to $\h$ as solution to \eqref{3}.
The first main result of this paper is the following theorem: 
\begin{thm}[Existence and uniqueness]\label{Thm1.1}
Let $R>0$.
There exists an open neighbourhood $\mathcal O$ of $0$ in $\V$ such that, for any initial data $\rho_0\in\mathcal O,$
there exists a maximal existence time $T:=T(\rho_0)>0$ and a unique classical solution  $\h=\rho(\cdot;\h_0)$ to problem \eqref{3} defined on $[0,T(\rho_0))$
which satisfies $\rho([0,T(\rho_0)))\subset \mathcal O$.
The mapping 
\[
\{(t,\rho_0)\,:\, \text{$\h_0\in\mathcal O$ and $0<t<T(\h_0)$}\}\mapsto\h(t;\h_0)\in h^{4+\alpha}(\s)
\]
is smooth.
\end{thm}

When $R=R_A$, we re-discover $\h\equiv0,$
 situation when the tumor is located at $D(0,R_A)$ as the unique radially symmetric stationary solution of \eqref{3}.
Concerning the stability properties of this solution we already know from the radially symmetric case \cite[Theorem 1.2]{EM} that this solution 
is unstable for $G<0$.
Moreover, we have:

\begin{thm}\label{P:G*G}
Let $R=R_A$ and $G_*>0$ be the constant defined by \eqref{eq:G*}.
Then the radially symmetric equilibrium $\h\equiv 0$ 
is unstable for all $G>G_*$. 
\end{thm}

Additionally to Theorem \ref{P:G*G} we have:
\begin{thm}\label{T:conv} Let $R=R_A$ and  assume that 
\begin{equation}\label{eq:ass}
\frac{A}{2} \frac{u_0'(1)}{u_0(1)}+A-f(1)> 0,
\end{equation}
whereby $u_0$ is the solution of \eqref{uniculu} for $n=0$.
Given $G>0,$ there exists a positive integer $l_G\in\N$ such that for all $\omega\in (0,\mu_0)$ and $l\geq l_G$ we find positive constants
$K_l>0$ and $\delta_l>0$ with the property that if $\|\h_0\|_{C^{4+\alpha}(\s)}\leq \delta_l$ and 
$\h_0$ is $2\pi/l-$periodic, then the solution $\h$ to \eqref{3} exists in the large, and
\[
\|\h(t)\|_{C^{4+\alpha}(\s)}+\|\h'(t)\|_{C^{1+\alpha}(\s)}\leq K_le^{-\omega t}\|\h_0\|_{C^{4+\alpha}(\s)},\quad t\geq0.
\]
Moreover, the solution $\h$ is $2\pi/l-$periodic for all $t\geq0.$
\end{thm}
We will show in the Appendix that the condition \eqref{eq:ass} is satisfied particularly when $f=\id_{[0,\infty)}$ and $R_A=1.$ 
 
 \section{The well-posedness result}
This section is dedicated to the proof of Theorem \ref{Thm1.1} and preparing Theorems \ref{P:G*G} and \ref{T:conv}.
A fundamental difficulty in treating problem \eqref{3}
is the fact that one has to work with unknown, variable domains $\0_\rho$.
We overcome this difficulty by transforming  problem \eqref{3} on the unitary disc $\0:= D(0,1)$.
Therefore, we define for all   $\rho\in\V$ the mapping
$\Theta_\rho:\R^2\to\R^2$ by
\[
\Theta_\rho(x)=Rx+\frac{Rx}{|x|}\f(|x|-1)\rho\left(\displaystyle\frac{x}{|x|}\right),
\]
where the cut-off function $\f\in C^\infty(\R,[0,1])$ satisfies
\[
\f(r)=\left\{
\begin{array}{llll}
&1,& |r|\leq 1/4,\\[2ex]
&0,& |r|\geq 3/4,
\end{array}
\right.
\]
and additionally $\max|\f'(r)|<4.$ 
Notice that $|x-1|\geq 3/4$  we have $\Theta_\rho(x)=x.$
Given  $x\in\s$, the mapping $[0,\infty)\ni r\mapsto r+\f(r-1)\rho(x/|x|)\in[0,\infty)$
is strictly increasing and therefore bijective.
The composition $\rho(\cdot/|\cdot|) $ has the same regularity properties as $\h$ on any subset of $\R^2$ which is bounded away from $0$, and using
the chain rule we have, cf. \cite{EM4}, that 
\begin{equation}\label{eq:grad2}
\nabla\left(\rho\left(\frac{x}{|x|}\right)\right)=\rho'\left(\frac{x}{|x|}\right)\left(-\frac{x_2}{|x|^2},\frac{x_1}{|x|^2}\right)
\end{equation}
for all $x\neq 0.$
Consequently, $\Theta_\rho $ is a diffeomorphism mapping $\0$ onto $ \0_\h$, i.e. $\Theta_\rho\in  \mbox{\it{Diff}}\,^{4+\alpha}(\0,\0_\rho)\cap \mbox{\it{Diff}}\,^{4+\alpha}(\R^2,\R^2).$
Such a diffeomorphism was first introduced by Hanzawa in \cite{EIH} to study the Stefan problem, and it is therefore called {\it Hanzawa diffeomorphism}. 
\begin{figure}\label{F:Hanz}
$$\includegraphics[width=0.7\linewidth]{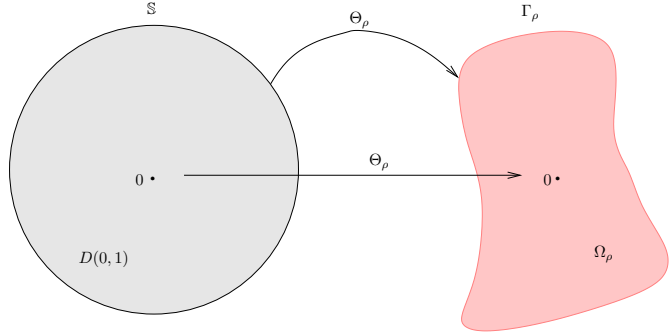}$$
\caption{The Hanzawa diffeomorphism}
\end{figure}
Additionally, we have that $\Theta_\rho(\s)=\Gamma_\rho$ (see Figure $2$).
The  push-forward operator induced by $\Theta_\rho$ is defined by
\[
\begin{array}{ll}
\Theta_\rho^*:\mbox{\it BUC}( \0_\rho)\to \mbox{\it BUC}(\0),& u\mapsto u\circ \Theta_\rho.
\end{array}
\]
These operators  allow us to transform the problem into an abstract Cauchy problem over $\s$.
General results of  the theory of maximal regularity, due to Sinestrari \cite{Si}, can be used to prove existence of a unique
classical solution, corresponding to small initial data.
The solution to \eqref{3} is then obtained (see Lemma \ref{Lem2.1} below) using the pull-back  operators defined by
 \[
\begin{array}{ll}
\Theta_* ^\rho:\mbox{\it BUC}(\0)\to \mbox{\it BUC}(\0_\rho),& v\mapsto v\circ \Psi_\rho,
\end{array}
\]
where $\Psi_\h:=\Theta_\h^{-1}=(\psi^1_\rho,\psi^2_\rho).$
The transformed  operators $\A(\rho)$ and $\B,$ are defined as follows.
Given $\h\in\V,$ $\A(\rho): buc^{2+\alpha}(\Omega)\to buc^\alpha(\Omega)$
is the differential operator given by 
\begin{equation}\label{eq:A(h)}
 \A(\rho):= \Theta_\rho^*\circ\Delta\circ \Theta_* ^\rho.
 \end{equation}
The operator $\A(\rho)$  is  linear and  uniformly elliptic, with  
\[
\A(\rho)v=b_{ij}(\h)v_{ij}+b_i(\h)v_i, \quad\forall v\in buc^{2+\alpha}(\0),
\]
whereby
\begin{align*}
b_{ij}(\rho)&=\psi^i_{\rho,1}(\Theta_\rho(x))\psi^j_{\rho,1}(\Theta_\rho(x))+\psi^i_{\rho,2}(\Theta_\rho(x))\psi^j_{\rho,2}(\Theta_\rho(x)),\\[1ex]
b_i(\rho)&=\psi^1_{\rho,11}(\Theta_\rho(x))+\psi^2_{\rho,22}(\Theta_\rho(x))
\end{align*}
for $ 1\leq i,j\leq2.$ 
Using \eqref{eq:grad2} and the chain rule, we can determine the coefficients $b_{ij}(\h)$ and $b_i(\h),$ $ 1\leq i,j\leq2,$
explicitly in terms of $\h$ and $\f,$ the cut-off function used when defining $\Theta_\h.$ 
Moreover, 
$\A$ depends analytically on $\h$
\begin{equation}\label{anl}
\A  \in C^\omega (\V, \kL( \mbox{\it buc}\, ^{2+\alpha}(\0), \mbox{\it buc}\, ^{\alpha}(\0))).
\end{equation}

The trace operator $ \B :\V \times buc^{2+\alpha}(\Omega)\times buc^{2+\alpha}(\Omega)\to h^{1+\alpha}(\s)$ is defined by the following relation
\begin{equation}\label{B(h)}
\B (\rho, v, q)= \frac{1}{R} \tr\left\langle G \nabla(\Theta_* ^\rho v)(\Theta_\rho)-\nabla(\Theta_* ^\rho q)(\Theta_\rho)
- \frac{AG}{2}\Theta_\rho,\  \nabla N_\rho(\Theta_\rho)\right\rangle,
\end{equation}
with $\tr$ the {\it trace operator} on $\s,$ i.e. $\tr v=v|_{\s}$ for $v\in \mbox{\it{BUC}}\,(\Omega),$
and  the curvature  $\kappa_{\Gamma_\rho}$ can be expressed in terms of $\h$  by the relation
\[
\kappa_{\Gamma_\rho}(\Theta_\h)=\frac{(1+\rho)^2+2\rho'^2-(1+\rho)\rho''}{R((1+\rho)^2+\rho'^2)^{3/2}}=:\kappa(\rho).
\]
It is not difficult to see that if $(\rho,\psi, p)$ is a solution of \eqref{3} then
$(\rho,v,q)=(\h,\Theta_\rho^* \psi,\Theta_\rho^* p)$ solves pointwise the following transformed problem
\begin{equation}\label{4}
\left \{
\begin{array}{rlllll}
\A(\rho)v&=&f(v)&\text{in}&\0,\\[0ex]
v&=&1 &\text{on}&\s,\\[0ex]
\A(\rho)q&=&0&\text{in}&\0,\\[0ex]
q&=& \kappa(\rho)-\displaystyle{\frac{AGR^2}{4}}(1+\rho)^2 &\text{on}&\s, \\[0ex]
\partial_t\rho&=&\B(\rho,v,q)&\text{on}&\s,\\[0ex]
\rho(0)&=&\rho_0.&&&
\end{array}
\right.
\end{equation}
The notion of solution for this problem is defined analogously to that of solution to \eqref{3}.
In fact the problems \eqref{3} and \eqref{4} are equivalent in the following sense:
 
\begin{lem}\label{Lem2.1} Given $\rho_0\in\V$ we have:\\[0.5ex]
$(a) $ If $(\rho,\psi, p)$ is a classical H\"older solution for \eqref{3}, then $(\rho,\Theta_\rho^*\psi,\Theta_\rho^*p)$ 
is a classical H\"older solution for \eqref{4}.\\[1ex]
$(b) $ If $(\rho,v,q)$ is a classical H\"older solution for \eqref{4}, then $(\rho,\Theta^\rho_*v,\Theta^\rho_*q)$ 
is a classical H\"older solution for \eqref{3}.
\end{lem}
\begin{proof}
The proof is  similar to the one   in \cite[Lemma 2.1]{EM4}. 
\end{proof}

\begin{lem}\label{Lem2.2}
The mapping
$$\V\ni\rho\mapsto \kappa(\rho)=\frac{(1+\rho)^2+2\rho'^2-(1+\rho)\rho''}{R((1+\rho)^2+\rho'^2)^{3/2}}\in h^{2+\alpha}(\s)$$
is analytic. 
Moreover, $\p\kappa(0)[\rho]=-(\rho+\rho'')/R,$ for all $\h\in h^{4+\alpha}(\s).$ 
\end{lem}
\begin{proof} The analyticity is obvious.
In order to compute the derivative one has only to calculate the gradient of a real valued function of three variables.  
\end{proof}
We introduce now solution operators to some semilinear, respectively linear Dirichlet problems related to our transformed problem \eqref{4}. 
From the  Leray-Schauder fixed point theorem (cf. \cite[Theorem 11.3]{GT01}) we obtain for each $\rho \in \V$ 
a  solution $u\in \mbox{ \it BUC}\,^{2+\alpha}(\0_\h)$ of problem 
\begin{equation}\label{qs}
\left \{
\begin{array}{rlllll}
\Delta u&=&f(u) &\text{in}&\0_\rho,\\[1ex]
u&=& 1 &\text{on}&\Gamma_\rho,\\[1ex]
\end{array}
\right.
\end{equation}
with $\rho \in \V.$
Using the maximum principle as we did in the proof of \cite[Theorem 2.6]{EM} we may prove the uniqueness of this solution.
Consequently, we have:
\begin{lem}\label{L.2.3}
Given $\rho\in\V,$ there exists a unique solution $\x(\rho)\in buc^{2+\alpha}(\0)$ of the semilinear Dirichlet problem
\begin{equation}\label{7}
\left \{
\begin{array}{rlllll}
\A(\rho)v&=&f(v) &\text{in}&\0,\\[1ex]
v&=& 1 &\text{on}&\s.\\[1ex]
\end{array}
\right.
\end{equation}
The mapping $\left[\V\ni\rho\mapsto \x(\rho)\in buc^{2+\alpha}(\0)\right]$
is smooth.
\end{lem}
\begin{proof}
For details we refer to the proof of \cite[Theorem 4.3.5]{AM}. 
\end{proof}

We consider now the solution operator corresponding to the second, linear Dirichlet problem in system \eqref{4}.
We state:
\begin{lem}\label{L.2.4}
Given $\rho\in\V,$ there exists a unique solution $\mS(\rho)\in buc^{2+\alpha}(\0)$ of the  Dirichlet problem
\begin{equation}\label{8}
\left \{
\begin{array}{rlllll}
\A(\rho)q&=&0 &\text{in}&\0,\\[1ex]
q&=& \kappa(\rho)-\displaystyle\frac{AGR^2}{4}(1+\rho)^2 &\text{on}&\s.\\[1ex]
\end{array}
\right.
\end{equation}
The mapping $[\V\ni\rho\mapsto \mS(\rho)\in buc^{2+\alpha}(\0)]$
is real analytic.
\end{lem}
\begin{proof}
Given $\h\in \V,$ the mapping
\[
(\A(\h), \tr): \mbox{ \it BUC}\,^{2+\alpha}(\0)\rightarrow \mbox{ \it BUC}\,^{\alpha}(\0)\times C^{2+\alpha}(\s)
\]
is a topological isomorphism from 
$\mbox{ \it BUC}\,^{2+\alpha}(\0)$ onto $\mbox{ \it BUC}\,^{\alpha}(\0)\times C^{2+\alpha}(\s).
$
It is well-known that the function mapping a bijective bounded linear operator onto its inverse is analytical; it can be expressed by a Neumann expansion in the neighbourhood of some other linear isomorphism. 

Hence, in view of Lemma \ref{Lem2.2} and equation \eqref{anl} it follows that 
\[
\mS(\h)=(\A(\h), \tr)^{-1} \left(0, \kappa(\rho)-\frac{AGR^2}{4}(1+\rho)^2\right)
\] 
is analytic. 
Since $\mS$ maps smooth functions on $\s$ into $\mbox{ \it BUC}\,^{\infty}(\0)$ we also have  $\mS(\h)\in \mbox{ \it buc}\,^{2+\alpha}(\0)$ for all $\h\in \V.$
\end{proof}
 
 \subsection{The nonlinear Cauchy problem}
 We use now the solution operators defined in Lemmas \ref{L.2.3} and  \ref{L.2.4} to transform the system \eqref{4} into an abstract Cauchy problem on the unit circle $\s$. 
We put  in the third equation of \eqref{4}  $\x(\rho)$, the solution to \eqref{7}, for $v$,  respectively  $\mS(\rho),$ the solution to \eqref{8}, for $q$,
to obtain  the following abstract Cauchy problem 
  \begin{equation}\label{9}
\begin{array}{ll}
\p_t\rho=\Phi(\rho),\qquad \rho(0)=\rho_0,
\end{array}
\end{equation}
where 
\begin{equation}\label{eq:Phi}
\Phi(\,\cdot\,):=\B(\,\cdot\,,\x(\,\cdot\,),\mS(\,\cdot\,))
\end{equation}
 is a nonlinear and nonlocal operator of third order which depends smoothly on $\rho$.
In order to prove Theorem \ref{Thm1.1} is suffices to show that  $\p\Phi(0)$ generates a strongly
continuous analytic semigroup in $\kL(h^{1+\alpha}(\s))$ with  definition domain
$h^{4+\alpha}(\s)$, that is 
\[
-\p\Phi(0)\in \kH(h^{4+\alpha}(\s),h^{1+\alpha}(\s)).
\]
The operator $\p\Phi(0)$ can be decomposed as the sum of its principal part, which has order three in $\h,$ with an operator of first order.
More exactly:
\begin{thm}\label{reguth}
The operator $\Phi$ is smooth, i.e. $\Phi \in C^\infty(\V, h^{1+\alpha}(\s)).$
Its derivative, $\p \Phi(0),$ writes  as the sum $\p\Phi(0)=A_1+A_2,$ where
\begin{equation}\label{a+b}
\text{$A_1\h:= \frac{1}{R^3}\p_\nu((\Delta, \tr)^{-1}(0,\h''))$ for $\h\in h^{4+\alpha}(\s),$}
\end{equation}
and $ A_2\in \kL(h^{2+\alpha}(\s),h^{1+\alpha}(\s)).$
\end{thm}

 In order to prove this theorem,  
we have to study first the regularity properties of the operator $\B,$ defined by \eqref{B(h)}.
It is convenient to we write this operator as a sum 
\[
\B(\rho,v,q)=\frac{G}{R} \B_1(\rho)v-\frac{1}{R} \B_1(\rho)q-\B_2(\rho),
\] 
where  $\B_1\in C^\omega(\V,\kL(\mbox{\it buc}^{2+\alpha}(\Omega), h^{1+\alpha}(\s))$ and $\B_2\in C^\omega(\V, h^{1+\alpha}(\s))$ are the operators defined by 
\begin{align*}
\text{$\B_1 (\rho) v= \tr\left\langle \nabla(\Theta_*^\rho v), \nabla N_\rho \right\rangle(\Theta_ \rho)$ and $\B_2 (\rho)= \frac{AG}{2R} \tr\left\langle \Theta_ \rho , \nabla N_\rho(\Theta_ \rho) \right\rangle.$}
\end{align*}
Since by the chain rule
\begin{align*}
\p\Phi(0)[\rho]=&-\frac{1}{R}\B_1(0)\p\mS(0)[\rho]-\frac{1}{R}\p\B_1(0)[\h]\mS(0)\\
&+\frac{G}{R}\p\B_1(0)[\h]\x(0)+\frac{G}{R}\B_1(0)\p\x(0)[\h]-\p\B_2(0)[\rho]
\end{align*}
for $ \h\in h^{4+\alpha}(\s),$ we must not only show that $\B_i,$ $i=1,2,$ have the regularity mentioned above but also determine their derivatives in $0$. 
We shall see that  the first term of this sum is the important one (corresponding to the operator $A_1$ in Theorem \ref{reguth})  since it is a third order operator,
and the last two terms are of lower order and play, as we shall see, no role when studying the well-posedness of the abstract evolution equation  \eqref{9}. 

Using relation \eqref{eq:grad2} we get  that
\begin{equation}\label{eq:gradd}
\nabla N_\h(\Theta_\h(x))=x-\frac{\h'(x)}{1+\h(x)}(-x_2,x_1),\quad x\in\s
\end{equation}
for all $\h\in\V.$
Particularly, we obtain find the following expression for $\B_2$
\begin{equation}\label{eq:B2grad}
\B_2(\h)=\frac{AG}{2R}\left\langle R(1+\h)x,x \right\rangle=\frac{AG}{2}(1+\h),
\end{equation}
wherefrom  we can easily see that $\B_2$ is analytic and that
\begin{equation}\label{13}
\p\B_2(0)[\rho]=\frac{AG}{2}\rho,\quad\forall \rho \in h^{4+\alpha}(\s).
\end{equation}
Consider now  the operator $\B_1$.
From the weak maximum principle we find that the function $\x(0)$ is radially symmetric, and one can easily see that $\mS(0)$ is constant.
Hence it suffices  to determine  $\p\B_1(0)[\h]v_0$  for a radially symmetric function $v_0\in \mbox{\it buc}^{2+\alpha}(\Omega) $ and $\h\in h^{4+\alpha}(\s).$

We state:
\begin{lem}\label{Lem3.1} The nonlinear operator $\B_1$ is analytic, i.e.  $$\B_1 \in C^\omega(\V\times buc^{2+\alpha}(\0),h^{1+\alpha}(\s)).$$
Given  $v_0\in buc^{2+\alpha}(\0)$ a radially symmetric function, we have that
\begin{equation}\label{10}
\p\B_1(0)[\cdot]v_0\equiv0.
\end{equation}
\end{lem}
\begin{proof}
Let $v_0\in buc^{2+\alpha}(\0)$ be a radially symmetric function and $\h\in\V.$
In view of \eqref{eq:gradd} we have that
\begin{align*}
\B_1(\h)& = \left[x_i\psi^j_{\h,i}(\Theta_\h)+\frac{\h'}{1+\h}
\left(x_2\psi^j_{\h,1}(\Theta_\h)-x_1\psi^j_{\h,2}(\Theta_\h)\right)
\right]\tr\p_j
\end{align*}
wherefrom we obtain the regularity assumption stated in the lemma.
Concering \eqref{10}  a detailed proof can be found in \cite[Lemma 4.4.2]{AM}.
\end{proof}

To finish the preparations for the proof of Theorem \ref{reguth} one more step must be done.
We have to determine the Fr\'echet derivative in $0$ of the analytic solution operator
defined in Lemma \ref{L.2.4}.

\begin{lem}\label{Lem3.4} 
Given $\rho\in\V$, the map $\p\mS(0)[\rho]\in\mbox{\it buc}^{2+\alpha}(\0)$ is the unique solution of the linear Dirichlet
problem
\begin{equation}\label{11}
\left \{
\begin{array}{rllll}
\Delta z &=& 0 &\text{in}&  \Omega,   \\[1ex]
z&=& -\left( \displaystyle \frac{1}{R}+\displaystyle\frac{AGR^2}{2}\right) \rho-\displaystyle\frac{1}{R}\rho'' &\text{on}&  \s.  
\end{array}
\right.
\end{equation}
\end{lem}
\begin{proof}
The proof is standard and we omit it.
\end{proof}

We come now to the proof of the main result of this subsection:
\begin{proof}[Proof of Theorem \ref{reguth}]
The regularity assumption follows directly from Lemmas \ref{L.2.3}, \ref{L.2.4},  \ref{Lem3.1}, and relation \eqref{eq:B2grad}.
Moreover since $\B_1(0)=R^{-1}\p_\nu$, we have that
\begin{align*}\label{dfi0}
\p \Phi(0)[\h]=&\frac{1}{R^3} \p_\nu\left((\Delta, \tr)^{-1}(0,\h'')\right)\\[1ex]
&+\left(\frac{1}{R^3}+\frac{AG}{2}\right)\p_\nu\left((\Delta, \tr)^{-1}(0,\h)\right)+\frac{G}{R^2}\p_\nu(\p \x(0)[\h])-\frac{AG}{2}\h.
\end{align*}
The  operator $\x$ defined  in Lemma \ref{L.2.3} can be extended to 
$\x:\{\h\in h^{2+\alpha}(\s)\,: \, \|\h\|_{C(\s)}<1/4\} \to \mbox{\it buc}^{2+\alpha}(\0),$ 
because we only need there that $\h$ is of class $C^{2+\alpha}$  to guarantee existence of a solution to \eqref{qs}.
Whence, the operator defined by 
\[
A_2\h:=\left(\frac{1}{R^3}+\frac{AG}{2}\right)\p_\nu\left((\Delta, \tr)^{-1}(0,\h)\right)+\frac{G}{R^2}\p_\nu(\p \x(0)[\h])-\frac{AG}{2}\h
\]
for $ \h\in h^{2+\alpha}(\s),$
belongs to $ \kL(h^{2+\alpha}(\s),h^{1+\alpha}(\s)).$ 
This completes the proof.
\end{proof}

We conclude this section with the proof of our first main result Theorem \ref{Thm1.1}.
\begin{proof}[Proof of the Theorem \ref{Thm1.1}]
The key point is showing that the operator $A_1$, which can be seen as the principal part of the Fr\'echet derivative $\p\Phi(0)$
generates a strongly continuous and analytic semigroup in   $\kL(h^{1+\alpha}(\s))$ for all $0<\alpha<1,$ i.e.
$-A_1\in \kH(h^{4+\alpha}(\s),h^{1+\alpha}(\s)).$
To this scope we show that $-A_1$ is a Fourier multiplier.
Given $\h\in h^{4+\alpha}(\s)$, we  consider its the Fourier expansion of $\rho=\sum_{k\in\Z}\wh\rho(k)x^k$, where $\wh \h(k):=\int_\s \h(x)x^k\, dx$ is the $k$-th Fourier coefficient of $\h.$
The well-known Poisson integral formula yields then 
$$(\Delta, \tr)^{-1}(0,\h'')=\sum_{k\in\Z} -k^2r^{|k|}\wh\rho(k)x^k$$
for all $r\leq 1$ and $x\in\s.$
Taking the derivative with respect to $r,$ in $r=1,$ we  finally obtain 
\begin{equation}\label{coeff}
A_1\left[\sum_{k\in\Z}\wh\rho(k)x^k\right]=\sum_{k\in\Z} \lambda_k\wh\rho(k)x^k
\end{equation}
for all $\sum_{k\in\Z}\wh\rho(k)x^k\in h^{4+\alpha}(\s),$ where 
or simplicity
\begin{equation}\label{a1}
\text{$\lambda_k:=\frac{-|k|^3}{R^3}$ for $k\in\Z$.}
\end{equation}
Following the same steps as in the proof of \cite[Theorem 3.5]{EM4} we obtain that 
$-A_1\in \kH(h^{4+\alpha}(\s),h^{1+\alpha}(\s))$ for all $\alpha\in(0,1).$

Since the constant $\alpha\in (0,1),$ fixed at the beginning of  Section 2 was arbitrary
we get that the assertions of Theorem \ref{reguth}  hold true with $\alpha$ replaced by $\beta,$ for some fixed  $\beta\in (0,\alpha).$
Particularly, we have that $-A_1\in \kH(h^{4+\beta}(\s),h^{1+\beta}(\s)),$
for some $\beta\in (0,\alpha).$ 
The definition domain of $A_2$ is an interpolation space between $ h^{1+\beta}(\s)$ and $h^{4+\beta}(\s),$ since
\[
h^{2+\beta}(\s)=(h^{1+\beta}(\s),h^{4+\beta}(\s))_{1/3}.
\]
We infer from \cite[Proposition 2.4.1]{L} that the sum $\p \Phi(0)=A_1+A_2$ also generates  a strongly continuous analytic 
semigroup in $\kL(h^{1+\alpha}(\s))$.
In view of \cite[Theorem 1.3.1]{Am}, the set $\kH(h^{4+\beta}(\s),h^{1+\beta}(\s))$ is open in
$\kL(h^{4+\beta}(\s),h^{1+\beta}(\s))$, whence there exists an open neighbourhood  
${\mathcal O}_\beta$ of $0$ in $h^{4+\beta}(\s)$ with the property that $-\p\Phi(\h)\in\kH(h^{4+\beta}(\s),h^{1+\beta}(\s))$
for all $\h\in {\mathcal O}_\beta.$
Then ${\mathcal O}:={\mathcal O}_\beta\cap h^{4+\alpha}(\s) $ is
 an open neighbourhood of $0$ in $\V.$
In view of  
\[
(h^{1+\beta}(\s),h^{4+\beta}(\s))_{\theta}=h^{1+\alpha}(\s).
\]
 we have establish that the assumptions of Theorem 8.4.1 in \cite{L}  hold
and the proof of Theorem \ref{Thm1.1} follows now from this theorem.
Consequently, given $\rho_0\in\mathcal O,$
there exists a positive time $T>0$ and a unique classical solution  $\rho$ to problem \eqref{3} on $[0,T]$
satisfying $\rho([0,T])\subset \mathcal O$.
Moreover, the solution may be extended on a maximal interval $[0,T(\rho_0))$ and if
$\rho$ is uniformly continuous with values in $h^{4+\alpha}(\s)$, then either
$$\underset{t\nearrow  T(\rho_0)}\lim \rho(t)\in\p{\mathcal O} \quad\text{or}\quad T(\rho_0)=+\infty.$$
This completes the proof of Theorem \ref{Thm1.1}.
\end{proof}

\section{Stability properties}

We study in this section the stability properties of the 
unique radially symmetric solution $D(0,R_A)$ determined in \cite{EM} when $A\in (0,f(1))$ and $G\neq0.$
Therefore, we choose  the constant $R$ fixed at the beginning of   Section 2 to be $R=R_A$.
Particularly, functions in $\V$ parametrise  domains near the stationary tumor $D(0,R_A)$.
We rediscover  $\h\equiv 0,$ situation when the tumor domain is the discus $D(0,R_A),$ as the unique radially symmetric stationary solution of \eqref{3}.

In order to study the stability properties of this equilibrium we have to determine the spectrum of the complexification of the  Fr\'echet derivative $\p\Phi(0),$
which we denote again by  $\p\Phi(0).$
The stability results established in Theorem \ref{P:G*G} and Theorem \ref{T:conv}
 are then obtained by applying the principle of linearised stability to  problem \eqref{9}.
 Repeating the arguments presented in the proofs of Theorem \ref{Thm1.1}  we see that the complexification of
 $\p\Phi(0)$ generates a strongly continuous and analytic semigroup.
Taking into consideration that the embedding $h^{4+\alpha}(\s,\C)\hookrightarrow h^{1+\alpha}(\s,\C)$ 
is compact, we deduce  that the complexification of
 $\p\Phi(0)$ has a compact resolvent.
From \cite[Theorem III.8.29]{K}, we conclude that its spectrum  consists only of eigenvalues of finite multiplicity,
\[
\sigma(\p\Phi(0))=\sigma_p(\p\Phi(0)).
\]
Given $\h\in h^{4+\alpha}(\s),$ we look in the following for the Fourier expansion of $\p\Phi(0)[\h].$
Having shown that $\p\Phi(0)$ is a Fourier multiplier,  the point spectrum of $\p\Phi(0)$ is given by the symbol of this multiplier.
The cornerstone of the  analysis leading to Theorems \ref{P:G*G} and \ref{T:conv} is the Theorem \ref{T:pmult}, which states that $\p\Phi(0)$
is a Fourier multiplier operator  with symbol explicitly determined. 

\begin{thm}\label{T:pmult} Given $\h\in h^{4+\alpha}(\s),$  we let $\h=\sum_{k\in\Z}\wh\h(k)x^k$ denote its associated Fourier series. 
 We have that 
\begin{equation}\label{eq:PHI}
\p\Phi(0)\left[\sum_{k\in\Z}\wh\h(k)x^k\right] 
=\underset{k\in \Z}\sum \mu_k \widehat \h(k) x^k, 
\end{equation}
where the symbol $(\mu_k)_{k\in\Z}$ is given by the relation
\begin{equation}\label{eq:symbol}
\mu_k:= -\frac{1}{R_A^3}|k|^3+ \frac{1}{R_A^3}|k| - G \left(\frac{A}{2}\frac{u_{|k|}'(1)}{u_{|k|}(1)}+A-f(1)\right),
\end{equation}
and $u_{|k|}$ is the solution of \eqref{uniculu} for $ n=|k|.$ 
Moreover,
\begin{equation}\label{eq:spectr}
\sigma(\p\Phi(0))=\{\mu_k\,:\, k\in\Z\}.
\end{equation}
\end{thm}

In order to obtain a Fourier expansion for the derivative $\p\Phi(0)[\h],$ $\h\in \V,$ we are left, cf. Theorem \ref{reguth} and \eqref{coeff}, to 
determine  the Fr\'echet derivative in $0$ of the  solution  operator $\x$ defined in Lemma \ref{L.2.3}.
However, the computations are more involved when computing $\p\x(0),$ since
$\x(0)$ is not a constant function.
We have that:

\begin{lem}\label{Lemma3.5}
Given $\rho\in h^{4+\alpha}(\s)$, the function $\p\x(0)[\rho]$ is the unique solution of the linear Dirichlet
problem
\begin{equation}\label{dS}
\begin{aligned}[b]
\Delta w&-R_A^2f'(v_0)w=\\
&-v_{0,11}  \frac{- 2x_2^2\f\rho-2x_1^2 |x| \f' \rho+2x_1 x_2\f\rho'}{|x|^3}\\
&- 2v_{0,12}\frac{2x_1x_2\f\h-2x_1x_2|x|\f'\h+(x_2^2-x_1^2)\f\h'}{|x|^3}\\
&+v_{0,22} \frac{2x_1^2\f\h+2x_2^2|x|\f'\h+2x_1x_2\f\h'}{|x|^3}\\
&-v_{0,1}\frac{x_1\f\h-x_1|x|\f'\h+2x_2\f\h'-x_1|x|^2\f''\h-x_1\f\h''}{|x|^3}\\
&-v_{0,2}\frac{x_2\f\h-x_2|x|\f'\h-2x_1\f\h'-x_2|x|^2\f''\h-x_2\f\h''}{|x|^3}  \,\quad {\rm in} \ \Omega,  \\
w &= 0 \quad {\rm on} \ \s, 
\end{aligned}
\end{equation}
where $v_0=\x(0)$ and $\f$ the cut-off function used to define the Hanzawa diffeomorphism.
\end{lem} 
\begin{proof}
The proof is standard though lengthy.
 For detailed  calculations  we refer to \cite[Lemma 5.1.1]{AM}.
\end{proof} 
The result of this lemma is not very useful yet.
This is due to the fact that the  first equation of \eqref{dS} 
contains besides $\h$ and $w$ also derivatives of  $\varphi$ and $v_0$.
Therefore, it is difficult to determine the Fourier expansion of $w,$ the solution of \eqref{dS}, when knowing
that of $\h.$ 
That is why we formally linearise the free boundary problem describing the stationary states of the full system \eqref{eq:problem} 
at the unique radially symmetric solution $(0, \psi_A, p_A)$, found in \cite[Theorem 1.1]{EM}, where we simply write $\psi_A:=\x(0)\circ\Theta_0$ 
for the solution of \eqref{qs} and   $p_A:=R_A^{-1}-AGR_A^2/4$ is the composition $p_A=\mS(0)\circ\Theta_0.$
By doing this we shall find a nice decomposition of the derivative $\p\x(0)$ as a sum of two operators (see Lemma \ref{expr} below).
Hence, we consider now perturbations of the radially symmetric solution of the form
\begin{equation*}
\left\{
\begin{array}{llll}
\psi _\varepsilon &=& \psi _A+\varepsilon \psi, \\[1ex]
p_\varepsilon &=& p_A+ \varepsilon p, \\[1ex]
\Omega _{\rho_\varepsilon}&=& \left\{re^{is }: 0\leq r< R_A(1+\varepsilon  \rho(s)), \, s\in\R \right\},
\end{array}
\right.
\end{equation*}
where we simply write $\h(s)=\h(e^{is})$ for all $s\in\R.$
Here, $\varepsilon$ is a small parameter, and $\psi, p, \rho $ are unknown functions.
The linearisation of problem \eqref{eq:problem} is then the following free boundary problem 
\begin{equation}\label{linpb}
\left \{
\begin{array}{rlllllll}
\Delta\psi &=& f'(\psi_A)\cdot \psi  &\text{in} \ D(0,R_A), \\[2ex]
\Delta p &=& 0 & \text{in}\ D(0,R_A), \\[2ex]
\psi &=& -R_A\psi_A'(R_A)\h\ & \text{on}\ \partial D(0,R_A), \\[2ex]
p&=& -\left(\displaystyle \frac{1}{R_A}+\displaystyle\frac{AGR_A^2}{2} \right)\rho-\displaystyle\frac{1}{R_A}\rho''&\text{on} \ \partial
D(0,R_A),\\[3.5ex]
G \p_\nu\psi- \p_\nu p &=& \displaystyle\frac{AGR_A}{2}\h
 & \text{on}\ \partial D(0,R_A).
\end{array}
\right.
\end{equation}
We look now for a connection between the 
linearisation \eqref{linpb} and the Fr\'echet derivative of $\Phi $  in $0$.
To this scope, we transform first the system  \eqref{linpb} to the unitary disc, i.e. we set 
$$w(x) =  \psi  ( R_Ax ), \,  z(x) = p  ( R_Ax ), \, v_0(x)= \psi_A( R_Ax),$$
for $x\in\0,$
and substitute these expressions  in \eqref{linpb}. 
This leads us to the following system of equations
\begin{equation}\label{linpbsubst}
\left \{
\begin{array}{rllll}
\Delta w &=& R_A^2f'(v_0)\cdot w  &\text{in}\quad  \0, \\[2ex]
w &=& -\partial_\nu v_0 \rho & \text{on}\quad\s, \\[2ex]
\Delta z &=& 0 & \text{in}\quad  \0,\\[2ex]
z &=& -\left(\displaystyle\frac{1}{R_A}+\displaystyle\frac{AGR_A^2}{2} \right)\rho-\displaystyle\frac{1}{R_A}\rho''&\text{on} \quad \s,\\[2.5ex]
\displaystyle\frac{G}{R_A} \partial_\nu\psi-\displaystyle\frac{1}{R_A} \partial_\nu p&=& \displaystyle\frac{AGR_A}{2}\rho 
 & \text{on}\quad \s.\\[2ex]
\end{array}
\right.
\end{equation}

Given $\rho \in h^{4+\alpha}(\s)$, we let $\W(\h) \in buc^ {2+\alpha}(\Omega)$ denote the solution to the linear Dirichlet problem
\begin{equation}\label{40}
\left \{
\begin{array}{rllll}
\Delta w &=& R_A^2f'(v_0)\cdot w &\text{in} \quad \0, \\[1ex]
w &=& -\partial_\nu v_0 \rho& \text{on}\quad \s.
\end{array}
\right.
\end{equation}
Further on, we want to determine a 
relation between $\p\x(0)[\h]$ and $\W(\h).$
 Therefore, we define the extension operator $\E:h^{4+\alpha}(\s)\rightarrow buc^{2+\alpha}(\Omega)$ by 
\begin{equation}\label{41}
\E(\h)(x):= \left( v_{0,1}(x)\frac{x_1}{|x|}+ v_{0,2}(x)\frac{x_2}{|x|}\right)\rho \left(\frac{x}{|x|}\right) \varphi(|x|-1)
\end{equation}
for $x\in\0.$
Using these operators we can now write  $\p\x(0)$ as the sum of $\W$ and $\E.$
This decomposition is very useful because we got rid, in this way, of all the terms from the right hand side of first equation of \eqref{dS}.
Recall that our goal is to determine the Fourier expansion of  $\p_\nu(\p\x(0)[\h])$ when $\h\in h^{4+\alpha}(\s).$
It turns out, cf. Lemma \ref{expr}, that $\p_\nu(\E(\h))$ is collinear with $\h$ for all $\h\in h^{4+\alpha}(\s).$
Moreover, using ODE--techniques we are able  to determine an expansion for
$\W(\h)$ for all $\h\in h^{4+\alpha}(\s), $ cf. \eqref{FourierW}.
Indeed, we have:
\begin{lem}\label{expr}
It holds that $\p\x(0)= \W+ \E.$
Moreover, given $\h\in h^{4+\alpha}(\s),$ we have that
\begin{equation}\label{A1}
\p_\nu(\E(\h))=\alpha_{A}\rho,
\end{equation}
where $\alpha_{A}:=\p^2U/\p r^2(1,R_A^2)>0$ and $U$  is the solution of  the parameter-dependent problem:
\begin{equation}\label{eq:param}
\left \{
\begin{array}{rlllll}
\displaystyle\frac{\p^2U}{\p r^2}(r,\lambda)+ \displaystyle\frac{1}{r}\displaystyle\frac{\p U}{\p r}(r,\lambda
 ) &=& \lambda f(U(r,\lambda )), & \ 0\leq r \leq 1, \\[2ex]
\displaystyle\frac{\p U}{\p r}(0,\lambda )& =& 0,\\[2ex]
U(1,\lambda )& =& 1,
\end{array}
\right.
\end{equation}
with $\lambda\in[0, \infty).$
\end{lem}
\begin{proof}
The proof follows by direct computation by taking into consideration that $v_0=\x(0)=U(|\cdot|, R_A^2).$
\end{proof}
 
In virtue of Lemma \ref{expr},  if we determine a Fourier expansion for   $\W(\h)$, then we completed the task of determining the expansion of  $\p\x(0)[\h]$ for all  $\h\in h^{4+\alpha}(\s).$ 
For this reasoning,  we consider expansions of the form
\begin{equation}\label{coeffi}
\begin{array}{rll}
\W(rx) &=& \underset{k\in \Z}\sum A_k(r)x^k\\[1ex]
\rho(x) &=& \underset{k\in \Z}\sum \widehat\rho(k) x^k
\end{array}
\end{equation}
with $r\in [0,1]$ and $x\in\s.$
 Substituting these expressions  into  \eqref{40}, and comparing the coefficients of $x^k$, we come 
to the following problem for the unknown function $A_k:$
\begin{equation}\label{iktheta}
\left \{ 
\begin{array}{rlll}
 A_k''+\displaystyle\frac{1}{r}A_k'-\displaystyle\frac{k^2}{r^2}A_k&=&R_A^2f'(v_0)A_k, \quad 0<r<1, \\[2ex]
 A_k(1)&=&- v_0'(1) \widehat\rho(k).
\end{array}
\right.
\end{equation} 
We have used here the relation $\p_\nu v_0 = v_0'(1)$ on $\s,$ where we identified $v_0$ with its restriction to the interval $[0,1].$
In order to prove the existence and uniqueness of the solution to \eqref{iktheta} we consider first the  associated problem \eqref{uniculu}.
The solution of \eqref{iktheta} will be then expressed in terms of the solution of this new system.
Thus, given $n\in \N$,  the problem 
\begin{equation}\label{uniculu}
\left\{
\begin{array}{rlll}
u''+\displaystyle\frac{2n+1}{r}u' &=& R_A^2f'(v_0)u, \quad & 0<r<1,\\[2ex]
u(0) &=& 1, \\[2ex]
u'(0)&=&0,
\end{array}
\right.
\end{equation}
has a unique solution $u_n \in C^\infty([0,1]).$\pagebreak

With this notation we have:
\begin{lem}\label{Ak}
Given $k\in \Z$, problem \eqref{iktheta} possesses a unique solution  $A_k \in C^\infty ([0,1])$ explicitly given by
\begin{equation}\label{Akrepr}
A_k(r)=\frac{-v_0'(1)}{u_{|k|}(1)}\widehat\rho(k)r^{|k|}u_{|k|}(r), \quad 0\leq r\leq1.
\end{equation}
We denoted by  $u_{n}, n\in\N,$   the  solution of \eqref {uniculu}.
\end{lem}
\begin{proof}
That $A_k$ is a solution of \eqref{iktheta} follows by direct computation.
The uniqueness may be obtain by a contradiction argument.  
\end{proof}

We give now a short proof of the theorem  stated at the beginning of this section.
\begin{proof}[Proof of Theorem \ref{T:pmult}]

From Lemma \ref{Ak} we obtain the following expansion for $\W(\rho)$
\begin{equation}\label{FourierW}
\W(\h)(rx^k)= \underset{k\in \Z}\sum \left[ -\frac{v_0'(1)}{u_{|k|}(1)}r^{|k|}u_{|k|}(r) \right ]\widehat\h(k) x^k
\end{equation}
for all $\h\in h^{4+\alpha}(\s).$
Let us determine now the constant $\alpha_A$ form Lemma \ref{expr}.
 From the first equation of   \eqref{eq:param}  we find that the constant $\alpha_A$ satisfies the relation $\alpha_{A}= R_A^2 f(1)- v_0'(1).$
Moreover, from the same equation and  relation (3.27) in\cite{EM} we have that 
\begin{equation}\label{eq:v0'}
v_0'(1)=AR_A^2/2,
\end{equation} 
hence
\begin{equation}\label{alfaR}
\alpha_A=R_A^2\left(f(1)-\frac{A}{2}\right).
\end{equation}
In view of Theorem \ref{reguth}, \ref{coeff}, and \eqref{FourierW} we conclude \eqref{eq:PHI} and the proof is complete.
\end{proof}

\subsection{Estimates for the symbol of the  derivative of $\bf \Phi$}

In order to study the stability properties of the unique radially symmetric equilibrium determined in \cite{EM} we need to study the sign of symbol 
$(\mu_k)_{k\in\Z}$ given by \eqref{eq:symbol} in dependence of the parameters $(A,G)\in(0,f(1))\times(0,\infty).$
We consider in here just the case when $G> 0, $ since for $G<0$ we already established in \cite[Theorem 1.2. (d)]{EM} 
that this circular equilibrium is 
unstable.
It is worth noticing that $\mu_k=\mu_{-k}$ for all $k\in\Z,$ so that we consider just the terms $\mu_k$ with $k\in\N.$
At first we have to as certain that $0$ is in  the spectrum of $\p\Phi(0)$
for all $(A,G)\in(0,f(1))\times[0,\infty).$

\begin{prop}\label{k1}
We have  that 
\[
\mu_1=0.
\]
\end{prop}
\begin{proof}
Since $v_0=\x(0)$ satisfies $v_0(r)=U(r,R_A^2)$ for all $0\leq r\leq 1,$
we have that
\[
v_0''+ \displaystyle \frac{1}{r} v_0' = R_A^2 f(v_0)\quad\text{in} \, (0,1). 
\]
Differentiating this equation with respect to $r$ and setting  $w_0:=v_0'$ we get
\begin{equation}\label{wo}
\left \{
\begin{array}{rlll}
w_0''+ \displaystyle\frac{1}{r}w_0'- \displaystyle\frac{1}{r^2} w_0 &=& R_A^2 f'(v_0)w_0,  \\[1ex]
w_0(1)& =& v_0'(1).
\end{array}
\right.
\end{equation}
Given $r\in[0,1]$, define $w(r)=v_0'(1) r u_1(r)/u_1(1),$
where $u_1$ denotes the solution of the problem \eqref{uniculu} when $n=1.$
With  $c:=v_0'(1) / u_1(1)$, we obtain by differentiation that
\[
w'(r) = c(u_1(r)+ ru_1'(r)),\qquad w''(r) = c(r u_1''(r)+2 u_1'(r)), 
\]
which in turn  implies that $w$ is a solution of the equation
\begin{align*}
w''+ \displaystyle\frac{1}{r}w'-\frac{1}{r^2}w= c(r u_1''(r)+3 u_1'(r))=R_A^2f'(v_0)w.
\end{align*}
Moreover, for $r=1$, we get  
$w(1)= v_0'(1).$
Since the solution of \eqref{wo} is unique,  we get $w_0=w,$
meaning that $v_0'(r)=v_0'(1) r u_1(r)/u_1(1)$ for all $0\leq r\leq 1.$
Rearranging this relation it follows that
\[
\frac{v_0'(r)}{v_0'(1)r}= \frac{u_1(r)}{u_1(1)},
\]
which leads to
\begin{align*}
\frac{u_1'(r)}{u_1(1)}&= \frac{v_0''(r)\, r \, v_0'(1)-v_0'(r)\, v_0'(1)}{v_0'(1)^2 \, r^2}= \frac{1}{v_0'(1)} \left[\frac{v_0''(r)}{r}- \frac{v_0'(r)}{r^2}\right ].
\end{align*}
If we choose $r=1$ this last relation leads to
\begin{align*}
\frac{u_1'(1)}{u_1(1)}&=  \frac{1}{v_0'(1)} \left[ v_0''(1)- v_0'(1) \right]= \frac{v_0''(1)}{v_0'(1)}-1.
\end{align*}
In view of relations \eqref{eq:v0'} and \eqref{alfaR}, and the definition of $\alpha_A,$ we have
\begin{align*}
\displaystyle\frac{u_1'(1)}{u_1(1)}&= \frac{-\displaystyle\frac{AR_A^2}{2}+R_A^2f(1)}{\displaystyle\frac{AR_A^2}{2}}-1= -2+\frac{2f(1)}{A}.
\end{align*}
Inserting this result in the expression of $\mu_1$ yields 
\[
\mu_1= -\frac{1}{R_A^3}+\frac{1}{R_A^3}-G \left[\frac{A}{2} \left(-2+\frac{2f(1)}{A}\right)+A-f(1)\right]=0.
\]
\end{proof}

Proposition \ref{k1} reveals that $\mu_1=0$ belongs to the spectrum of $\p\Phi(0)$.
However, we show now that the sequence $\mu_k\to_{k\to\infty}-\infty.$
This is obviously true if the sequence $(u_n'(1)/u_n(1))_{n\in\N}$ is bounded. 
Even more, we have:
\begin{lem}\label{un+1}
It holds that 
\[
\displaystyle\frac{u_n'(1)}{u_n(1)} \underset{n \to \infty}\longrightarrow 0,
\]
\end{lem}
\begin{proof}
Let $n\in \N$ and set $v:= u_{n+1}-u_n \in C^\infty ([0,1]).$ 
Recall that  $u_n$ is an increasing  function for all $n\in\N.$
This can be seen  since
\begin{align}\label{50}
u_n(r)=1+ \int_0^r \frac{R_A^2}{s^{2n+1}}\int_0^s \tau^{2n+1}f'(v_0(\tau))u_n(\tau)\, d\tau \, ds, \quad 0\leq r\leq 1,
\end{align}
relation which is obtained by multiplying the first equation of problem \eqref{uniculu} with $r^{2n+1}$ and then integrating twice.
From  \eqref{uniculu} we obtain
\begin{equation*}
\left \{
\begin{array}{rlll}
v''+\displaystyle\frac{2n+1}{r}v'&=& R^2 f'(v_0)v-\displaystyle\frac{2}{r}u_{n+1}', &0< r< 1,\\[1ex]
v'(0)&=& 0,\\[1ex]
v(0)&=& 0.
\end{array}
\right.
\end{equation*}
Furthermore,  we have
\begin{align*}
\underset{t\to 0}\lim \frac{v'(t)}{t} &= \underset{t\to 0}\lim \frac{u_{n+1}'(t)-u_n'(t)}{t}= R_A^2f'(v_0(0))\left( \frac{1}{2n+4}-\frac{1}{2n+2}\right)\\[1ex]
&= -R_A^2 f'(v_0(0)) \frac{2}{(2n+4)(2n+2)} < 0,
\end{align*}
which implies by \eqref{eq:conditions} that $v'(t) <0$ for $t\in (0, \delta)$ and  some  $\delta<1.$ 
Thus, $v$ is decreasing on $(0, \delta).$
Let now $t\in [0,1]$  and $m_t:= \min_{[0,t]} v \leq 0.$ 
A maximum principle argument shows   that  the non-positive minimum must be achieved at $t$,  $m_t=v(t)$ which implies $u_{n+1}(t)\leq u_n(t)$ for all $t\in [0,1]$.
Particularly, $u_{n+1}(1) \leq u_n(1)$.
It also holds that
\begin{align*}
u_n'(t)&= R_A^2 \int_0^t \left(\frac{\tau}{t}\right)^{2n+1}f'(v_0(\tau))u_n(\tau) \, d\tau \\[1ex]
& \geq  R_A^2 \int_0^t \left(\frac{\tau}{t}\right)^{2n+3}f'(v_0(\tau))u_{n+1}(\tau) \, d\tau \\[1ex]
&= u_{n+1}'(t)
\end{align*}
for all $t\in [0,1],$
hence $
u_{n+1}'(1)\leq u_n'(1)$.
Moreover, from
\[
u_n'(t)=\frac{R_A^2 \displaystyle\int_0^t \tau^{2n+1}f'(v_0(\tau))u_n(\tau) \, d\tau }{t^{2n+1}}
\]
it follows that
\[
u_n'(1)= R_A^2\int_0^1 \tau^{2n+1}f'(v_0(\tau))u_n(\tau) \, d\tau.  
\]
The dominated convergence theorem implies $u_n'(1)\searrow 0.$
 A similar argument provides  $u_n(1)\searrow 1.$ 
Therefore, $u_n'(1)/u_n(1)\rightarrow 0,$ and we are done.
\end{proof}

\subsection{Instability for $\bf G>G_{*}$}
We prove now that if $G$ is large enough, then the radially symmetric equilibrium $D(0,R_A)$, determined in  \cite[Theorem 1.1]{EM},
is unstable, in the sense that problem \eqref{9} has  nontrivial backwards solutions.
This reveals that the exponential stability result stated in \cite[Theorem 1.2. (a)]{EM} gives a false impression 
about the stability properties of the radially symmetric equilibrium.
If $k\in\N$ is large enough, we find in view of the Lemma \ref{un+1} a unique value $G_k$ of the parameter $G$ such that $\mu_k=0$ iff $G=G_k.$
More precisely, given $k\in\N$ such that
\[
\frac{A}{2} \frac{u_k'(1)}{u_k(1)}+ A-f(1) \neq 0,
\]
we set
\begin{align}\label{Gk}
G_k:=\displaystyle\frac{\displaystyle{-\frac{1}{R^3}k^3+\frac{1}{R^3}k}}{\displaystyle{\frac{A}{2} \frac{u_k'(1)}{u_k(1)}+ A-f(1)}}.
\end{align} 
To prove the statement of Theorem \ref{P:G*G}, we set 
\begin{equation}\label{eq:G*}
G_*:=\min\left\{ G_k\,:\,  \frac{A}{2} \frac{u_k'(1)}{u_k(1)}+ A-f(1) <0\right\}\geq 0.
\end{equation}
Moreover, since $G_k\to_{k\to\infty}\infty$, the minimum must be achieved, i.e.
we find at least an integer $k_0\in\N$ such that 
\begin{equation}\label{eq:cond<}
\frac{A}{2} \frac{u_{k_0}'(1)}{u_{k_0}(1)}+ A-f(1)<0,
\end{equation}
and $G_*=G_{k_0}.$
With these preparations we are able to prove the instability result.
\begin{proof}[Proof of Theorem \ref{P:G*G}]
Let  $G>G_{*}$ be given and $k_0$ be an integer such that $G_*=G_{k_0}$  and \eqref{eq:cond<} holds true.
It follows that 
\begin{align*}
-\frac{1}{R^3} k_0^3+\frac{1}{R^3} k_0> G \left(\frac{A}{2} \frac{u_{ k_0}'(1)}{u_{ k_0}(1)}+ A-f(1) \right),
\end{align*}
which implies $\mu_{k_0}>0.$
We are left to check the following  instability assumptions 
\begin{equation*}
\left\{
\begin{array}{rlll}
\sigma_+(\p \Phi (0))= \sigma (\p \Phi (0)) \cap \{\lambda \in \C : \Ro \lambda >0\} &\not =& \emptyset,\\[1ex]
\inf\{\Ro \lambda : \lambda \in \sigma_+(\p \Phi (0))\}&>&0,
\end{array}
\right.
\end{equation*}
where $\p\Phi(0)$ stands here for the complexification of $\p\Phi(0)$.
The first one is clear  since  $\mu_{ k_0}>0$.
Moreover, we infer from  Lemma \ref{un+1} that  $\mu_{k}\rightarrow_{k\to\infty} -\infty,$ and therefore the unstable spectrum  $\sigma_+(\p \Phi (0))$ 
contains finitely many   positive eigenvalues of the Fr\'echet derivative $\p\Phi(0)$. 
Thus, we found out that the assumptions of \cite[Theorem 9.1.3]{L}  are satisfied and therewith the radially symmetric solution $D(0,R_A)$ is unstable.
\end{proof}

\subsection{Exponential convergence of $\bf 2\pi/l-$periodic data}
In this section we show that the exponential stability result stated in \cite[Theorem 3.0.3 (a)]{AM} can be generalised for
 solutions of the original problem \eqref{eq:problem} which correspond to initial data $\h_0$  closed to the unique radially symmetric equilibrium 
and $2\pi/l$--periodic. 
The positive integer $l$ depends on the constant $G,$ which is chosen to be positive, so that $D(0,R_A)$
is the unique equilibrium of the problem \eqref{eq:problem}. 
The main result of this section, Theorem \ref{T:conv} requires the following assumption
\begin{equation*}
\frac{A}{2} \frac{u_0'(1)}{u_0(1)}+A-f(1)> 0,
\end{equation*}
which means that the eigenvalue $\mu_0=\mu_0(G)$ is negative for all $G>0.$
This assumption is satisfied for example if $R_A=1$ and $f=\id_{[0,\infty)}$ (see Appendix).

In order to prove  Theorem \ref{T:conv} we introduce first appropriate subspaces of the small H\"older spaces.
 Given $k\in\N$ and $l\in\N, n\geq2$, we define the  subspace of $h^{k+\alpha}(\s)$ consisting of $2\pi/l$--periodic  functions by
 \[
h^{k+\alpha}_{l}(\s):=\{\h\in h^{k+\alpha}(\s)\,:\, \text{$\h(x)=\h(e^{2\pi i/l}x)$ for all $ \, x\in\s$}\}.
\]
Set further $\V_{l}:=\V\cap h^{4+\alpha}_{l}(\s).$
The Fourier series associated to  $\h\in h^{k+\alpha}_{l}(\s)$ is 
\[
\h=\sum_{k=0}^\infty \wh\h(kl) x^{kl},
\]
where $\wh\h(kl)$ is the $kl-$th Fourier coefficient of $\h.$
Our  first objective  is to prove that if $\h\in  \V_{l},$ then $\Phi(\h)\in h^{1+\alpha}_{l}(\s),$ that is 
\begin{equation}\label{eq:rst1}
\Phi\in C^\infty(\V_l, h^{1+\alpha}_l(\s)).
\end{equation}
Having shown \eqref{eq:rst1}, by choosing $l$ large enough we can exclude the eigenvalues $\mu_k$ with $k$ small from the spectrum of $\p\Phi(0)$.
These are the eigenvalues which we could not estimate whether they are negative or not.
In this way we also eliminate  $\mu_1=0$ from the spectrum.

Let $\h\in\V_l$ be given and let  $\psi:=\Theta^\rho_*v$, where $v$ is the solution  of \eqref{7}.
We prove  that $\psi$ satisfies $ \psi(x)=\psi(e^{2\pi i/l}x)$ for all $x\in\0_\h.$ 
Therefore, we must prove first that, if $x \in \Omega_{\rho}$,  then $xe^{2\pi i/l}$ belongs to $\Omega_\rho.$
Indeed, given $x\in\0_\h$, we have that
\begin{align*}
|xe^{2\pi i/l}|=|x|< R\left(1+\rho\left( \frac{x}{|x|}\right)\right) = 
 R\left(1+\rho\left( \frac{xe^{2\pi i/l}}{|xe^{2\pi i/l}|}\right)\right),
\end{align*}
which implies that $xe^{2\pi i/l}\in \Omega_\rho.$ 
Recall that the function  $\psi \in \mbox{\it buc}^{2+\alpha}(\Omega_\h)$ is the unique solution of the Dirichlet problem  
\begin{equation}\label{Dirfi} 
\left \{
\begin{array}{rlcl}
\Delta \psi &=& f(\psi )  &\text{in} \ \Omega _{\rho},  \\[1ex]
\psi &=& 1 & \text{on}\ \Gamma _{\rho}.
\end{array}
\right.
\end{equation}
Setting $\ov \psi (x):= \psi (xe^{2\pi i/l})$ for $ x\in\0_\h,$ 
we have defined in this way a further  solution of \eqref{Dirfi}, since 
\[
\Delta \ov \psi(x) = \Delta \psi (xe^{2\pi i/l})= f(\psi (xe^{2\pi i/l}))= f(\ov \psi(x)),
\]
for all $x\in\0_\h$ and $\ov \psi=1$ on $\Gamma_{\h}$.
The uniqueness of the solution to \eqref{Dirfi}   implies  that $\psi = \ov \psi.$
 Thus 
$\psi (x)= \psi (xe^{2\pi i/l})$
for all $x\in \Omega_\h.$
Following the same schema we can prove that $p:=\Theta^\rho_*q$,  where $q$ is the solution  of \eqref{8}, satisfies 
$ p(x)=p(e^{2\pi i/l}x)$ for all $x\in\0_\h$ provided $\h\in\V_{l}.$
With these preparations we state:

\begin{lem}\label{L:rest}
Given $l\geq 2$, the operator $\Phi$ maps smoothly $ \V_{l}$ into $h^{1+\alpha}_{l}(\s).$
\end{lem}
\begin{proof}
It remains to show that $\B_1(\h) \x(\h)$ and $\B_1(\h) \mS(\h)$ are $2\pi/l$--periodic.
The assertion follows then in view of \eqref{eq:Phi} and relation \eqref{eq:B2grad}.
 We prove just the assertion  for $\B_1(\h)\x(\h)$, the proof that  $\B_1(\h, \mS(\h))$ is  $2\pi/l$--periodic  follows analogously. 
Indeed, given $x\in \s$ we have
\begin{align*}
\B_1(\h, \x(\h))(x) =& \langle \nabla (\Theta^\h_* v)(\Theta_{\rho}(x)) , \nabla N_\h(\Theta_\h (x))\rangle \\[1ex]
=& \left\langle \nabla \psi(\Theta_\h(x)), \frac{\nabla N_\h(\theta_\h (x))}{|\nabla N_\h(\Theta_\h (x))|}\right\rangle  
|\nabla N_\h(\Theta_\h (x))|\\[1ex]
=& \, \p_\nu \psi(\Theta_\h(x))  |\nabla N_\h(\Theta_\h (x))|,
\end{align*}
where $\psi$ is the solution of \eqref{Dirfi}.
Additionally, $\left|\nabla N_\h\left(\Theta_\h\left(xe^{2\pi i/l}\right)\right)\right|=|\nabla N_\h(\Theta_\h(x))|$
for all $x\in\s.$
In order to prove that also $\p_\nu \psi(\Theta_\h(xe^{2\pi i/l}))=\p_\nu \psi(\Theta_\h(x))$, $x\in\s,$  we introduce the rotation matrix 
\[
M=\left( \begin{array}{c}
\begin{array}{cc}
\cos\left(\frac{2\pi}{l}\right) & -\sin\left(\frac{2\pi}{l}\right) \\[2ex]
\sin\left(\frac{2\pi}{l}\right) & \cos\left(\frac{2\pi}{l}\right)
\end{array}
\end{array}
\right).
\]
We infer,  from $\psi(x)= \psi\left(x e^{2\pi i/l}\right)= \psi(M\cdot x),$ that $\nabla \psi(x)= M^\top\cdot \nabla \psi(M\cdot x)$
for all $x\in\0_\h.$
Thus,
\begin{align*}
\p_\nu \psi \left(\Theta_\h \left(x e^{2\pi i/l}\right) \right)
&= \left\langle \nabla \psi \left(\Theta_\h \left(M \cdot x\right) \right),\nu_\h 
\left(\Theta_\h \left(M \cdot x\right) \right)\right\rangle\\[1ex]
&=\left \langle \nabla \psi \left(M\cdot\Theta_\h(x)\right), M\cdot\frac{x - \h'(x)  (-x_2, x_1)}
{\sqrt{(1+\h)^2+\h'^2}(x)}\right\rangle\\[1ex]
&=\left\langle  M \cdot\nabla \psi(\Theta_\h(x)), M\cdot\frac{x - \h'(x)  (-x_2, x_1)}{\sqrt{(1+\h)^2+\h'^2}(x )} \right\rangle\\[1ex]
&=\left \langle \nabla \psi(\Theta_\h(x)), M^\top\cdot M\cdot\frac{x - \h'(x)  (-x_2, x_1)}{\sqrt{(1+\h)^2+\h'^2}(x )} \right\rangle\\[1ex]
&=\left \langle \nabla \psi \left(\Theta_\h \left(x \right) \right),\nu_\h \left(\Theta_\h \left(x \right) \right)\right\rangle\\[1ex]
&=\p_\nu \psi \left(\Theta_\h \left(x \right) \right)
\end{align*}
for all $x \in \s$. 
Summarising, $\B_1(\h, \x(\h))\left(xe^{2\pi i/l} \right) = \B_1(\h, \x(\h))\left(x \right) $ for all $x\in\s, $ and the proof is completed.\\[2ex]
\end{proof}
The Fr\'echet derivative $\p\Phi(0)$ of the mapping $\Phi\in C^\infty(\V_l, h^{1+\alpha}_l(\s))$ is, in view of \eqref{eq:PHI}, given by the relation
\[
\p\Phi(0)\left[\sum_{k\in\Z}\wh\h(kl)x^{kl}\right] 
=\underset{k\in \Z}\sum \mu_{kl} \widehat \h(kl) x^{kl}, 
\]
where $(\mu_{kl})_k$ are defined by \eqref{eq:symbol}.
We come now to the proof of the exponential stability result for $2\pi/l-$ periodic data:

\begin{proof}[Proof of Theorem \ref{T:conv}]
Let $G>0$ be given.
Since $\mu_{|k|}\to_{|k|\to\infty}-\infty$, we find a positive integer $l_G$ such that
$\mu_{|k|}\leq \mu_0$ for all $|k|\geq l_G.$
Let $l\geq l_G$ be fixed.
In view of relation Lemma  \ref{L:rest},
we find that the restriction $\Phi\in C^\infty(\V_l,h^{1+\alpha}_l(\s))$
satisfies the assumptions of \cite[Theorem 9.1.2]{L}.
Indeed, since $l\geq l_G,$ it holds that $\mu_{kl}\leq \mu_0$ for all $k\in\N.$
Consequently,  the spectrum of the  complexification of the $\p\Phi(0)$
consists  only of  the negative eigenvalues $\{\mu_{kl}\,:\, k\in\N\}$,
and
is bounded away from the positive half plane by $\mu_0.$
The assertion follows now immediately  from \cite[Theorem 9.1.1]{L}.
\end{proof}

\section{Appendix}
We show now that the condition \eqref{eq:ass}, meaning that $\mu_0(G)<0$ for all $G>0,$ is not to restrictive.
\begin{obs}\label{O:o}
The assertion
\[
\frac{A}{2} \frac{u_0'(1)}{u_0(1)}+A-f(1)> 0
\]
is fulfilled when $f= {\rm id}_{[0,\infty)}$ and $R_A=1.$
\end{obs}
\begin{proof}
In view of Proposition \ref{k1}, our assertion is equivalent with 
\[
\frac{A}{2} \frac{u_0'(1)}{u_0(1)}+A-f(1)> 0=\mu_1(G)=\frac{A}{2} \frac{u_1'(1)}{u_1(1)}+A-f(1).
\]
Consequently, we have to show only that
\[
\frac{u_0'(1)}{u_0(1)}> \frac{u_1'(1)}{u_1(1)}.
\]

We assume now that  $u_0$, the solution of \eqref{uniculu} when $n=0$, is  analytic and the Taylor series associated to $u_0$ in $0$ 
\[
u_0= \sum_{k=0}^\infty a_k x^k,
\]
converges on $[0,1].$
Problem \eqref{uniculu} writes now as follows
\begin{equation*}
\left \{
\begin{array}{rlllll}
xu_0''+u_0'-xu_0 &=&0,   & 0\leq x\leq1, \\[1ex]
u_0(0)&=&1, \\[1ex]
u_0'(0)&=&0.
\end{array}
\right.
\end{equation*}
From the initial conditions of \eqref{uniculu} it follows immediately that $a_0=u_0(0)=1$ and $a_1=u_0'(0)=0.$ 
Plugging  $u_0$ and its derivatives in the first equation of the system, one finds out that
\begin{align}\label{anlyu0}
\sum_{k=1}^\infty k(k+1)a_{k+1} x^k+ \sum_{k=0}^\infty (k+1)a_{k+1} x^k-a_0x-\sum_{k=2}^\infty a_{k-1} x^k=0.
\end{align}
Identifying the coefficient of $x^k$ in \eqref{anlyu0} yields 
\begin{align*}
& a_2=\frac{a_0}{4}=\frac{1}{4},\\[1ex]
& a_{k+1}=\frac{a_{k-1}}{(k+1)^2}, \, \forall k\in \N,
\end{align*} 
and, from  $a_1=0$, we deduce that 
\begin{align*}
& a_{2k+1}=0 \quad\text{and}\quad a_{2k}= \prod_{n=1}^k \frac{1}{(2n)^2}, \, \forall k\in \N,
\end{align*}
thus, 
\begin{equation}\label{un0}
u_0(x)=1+\sum_{k=1}^\infty \left(\prod_{n=1}^k \frac{1}{(2n)^2}\right)x^{2k},\quad x\in[0,1].
\end{equation}

We make now the same assumption on $u_1,$ the solution of \eqref{uniculu} when $n=1.$
We  then get, that $u_1$ is the solution of the following system
\begin{equation*}
\left \{
\begin{array}{rlllll}
xu_1''+3u_1'-xu_1 &=&0,&0\leq x\leq 1, \\[1ex]
u_1(0)&=&1, \\[1ex]
u_1'(0)&=&0.
\end{array}
\right.
\end{equation*}
As above, we obtain 
\begin{align*}
\sum_{k=1}^\infty k(k+1)a_{k+1} x^k+ 3\sum_{k=0}^\infty (k+1)a_{k+1} x^k-a_0x-\sum_{k=2}^\infty a_{k-1} x^k=0,
\end{align*}
and therefore
\begin{align*}
& a_{2k+1}=0 \quad\text{and}\quad a_{2k}= \prod_{n=1}^k \frac{1}{2n(2n+2)}, \, \forall k\in \N.
\end{align*}
Thus,
\begin{equation}\label{un1}
u_1(x)=1+\sum_{k=1}^\infty \left(\prod_{n=1}^k \frac{1}{2n(2n+2)}\right)x^{2k},\quad x\in[0,1].
\end{equation}
It is worth  noticing that the Taylor series  associated to $u_0$ and $u_1$, respectively, in $0$ define analytic functions on the whole real line,
 so that the representations \eqref{un0} and \eqref{un1} are valid.

With three exact decimals we have that 
\[
\frac{u_0'(1)}{u_0(1)}\approx 0.446>0.240\approx \frac{u_1'(1)}{u_1(1)},
\]
which leads to the desired conclusion.
\end{proof}

\vspace{1cm}

\end{document}